\documentclass[reqno]{amsart}
\usepackage[margin=3cm]{geometry}
\usepackage{amsmath,amssymb,amscd}
\usepackage{amsthm}
\usepackage{ae} % pekné Ľ
\usepackage{enumitem}
\usepackage{xcolor, graphicx}
\usepackage{subcaption}
\usepackage{tikz}
\usepackage[normalem]{ulem}
\usepackage{array}
\usepackage{tikz}
\usetikzlibrary{arrows}
\usepackage{pgfplots}

\newcommand{\diam}{\mathop{\rm diam}}

\newcommand{\id}{\mathrm{id}}

\newcommand{\dist}{\mathop{\rm dist}}

\theoremstyle{plain}
\newtheorem{theorem}{Theorem}
\newtheorem{lemma}[theorem]{Lemma}
\newtheorem{proposition}[theorem]{Proposition}
\newtheorem{corollary}[theorem]{Corollary}
\theoremstyle{definition}
\newtheorem{definition}[theorem]{Definition}
\newtheorem{example}[theorem]{Example}
\newtheorem{remark}[theorem]{Remark}

\newtheorem{claim}{Claim}[theorem]

\title{Spaces where all closed sets are $\alpha$-limit sets}

\author{Jana Hant\'{a}kov\'{a}}
\address{Faculty of Applied Mathematics, AGH University of Science and Technology, al. Mickiewicza 30, 30-059 Krak\'{o}w, Poland}
\address{Mathematical Institute of the Silesian University in Opava, Na Rybni\v{c}ku 1, 74601, Opava, Czech Republic}
\email{jana.hantakova@math.slu.cz}

\author{Samuel Roth}
\address{Mathematical Institute of the Silesian University in Opava, Na Rybni\v{c}ku 1, 74601, Opava, Czech Republic}
\email{samuel.roth@math.slu.cz}

\author[{\v L}. Snoha]{{\v L}ubom\'\i r Snoha}
\address{Department of Mathematics, Faculty of Natural Sciences, Matej
	Bel University, Tajovsk\'eho 40, 974 01 Bansk\'a Bystrica, Slovakia}
\email{lubomir.snoha@umb.sk}

\makeatletter
\@namedef{subjclassname@2020}{%
  \textup{2020} Mathematics Subject Classification}
\makeatother
\subjclass[2020]{Primary: 37B02, Secondary:  37B45, 37B05}
\keywords{$\alpha$-limit set, space with enough arcs, AF-space.}

\makeatletter
\def\blfootnote{\xdef\@thefnmark{}\@footnotetext}
\makeatother

\let\oldtocsection=\tocsection
\let\oldtocsubsection=\tocsubsection
\renewcommand{\tocsection}[2]{\hspace{0em}\oldtocsection{#1}{#2}}
\renewcommand{\tocsubsection}[2]{\hspace{2em}\oldtocsubsection{#1}{#2}}

\begin{document}

\begin{abstract} 
Metrizable spaces are studied in which every closed set is an $\alpha$-limit set for some continuous map and some point. It is shown that this property is enjoyed by every space containing sufficiently many arcs (formalized in the notion of a space with enough arcs), though such a space need not be arcwise connected. Further it is shown that this property is not preserved by topological sums, products and continuous images and quotients. However, positive results do hold for metrizable spaces obtained by those constructions from spaces with enough arcs.
\end{abstract}

\maketitle

\blfootnote{\copyright 2022. This is author accepted manuscript of article published by Elsevier in Topology and its Applications (2022), doi: 10.1016/j.topol.2022.108035. The manuscript version is made available under the CC-BY-NC-ND 4.0 license https://creativecommons.org/licenses/by-nc-nd/4.0/}
\section{Introduction}\label{S:intro}

Limit sets of forward trajectories, known as $\omega$-limit sets,  play an important role in topological dynamics. They are well-studied, even a complete topological characterization of them is known for continuous maps in some spaces, see~\cite{Sp} and references therein. For $\alpha$-limit sets the situation is different. Of course for homeomorphisms there is nothing new to say, since in this case the $\alpha$-limit sets are just the $\omega$-limit sets of the inverse map. But very little is known for non-invertible maps. This is perhaps related to the fact that $\alpha$-limit sets are rather flexible objects. Other kinds of limit sets, like $\omega$-limit sets and branch $\alpha$-limit sets are always internally chain transitive, and this leads to strong topological restrictions on what sets can arise this way \cite{Hirsch}. But $\alpha$-limit sets need not be internally chain transitive, or even contained in the chain-recurrent set of the system. The only obvious general facts are that $\alpha$-limit sets must be closed and forward invariant. It seems that the problem of characterizing $\alpha$-limit sets, even in simple topological spaces, has not yet been considered. There is only the simple observation made in \cite{KMS} that in the unit interval $[0,1]$, every closed set can occur as an $\alpha$-limit set. Thus our main goal is to answer the question: \emph{For which topological spaces $X$ can all closed subsets of $X$ occur as $\alpha$-limit sets?} Note that we allow the map to vary, requiring only that each closed subset of $X$ is an $\alpha$-limit set for \emph{some} continuous self-map $f:X\to X$. When this holds, we will say that $X$ is an \emph{$\mathcal{AF}$-space}, since the two classes $\mathcal{A}(X)$ and $\mathcal{F}(X)$ of $\alpha$-limit sets and closed sets in $X$ coincide.

It turns out that the interval is far from the only $\mathcal{AF}$-space. We start by showing that the $\mathcal{AF}$-property is enjoyed by every space with \emph{enough arcs}, by which we mean a non-degenerate metrizable space $X$ in which each non-degenerate proper closed subset can be joined to its complement by an arc in $X$. Any space with enough arcs is necessarily connected. On the other hand, in the class of spaces which have enough arcs (and so are $\mathcal{AF}$-spaces) there are, trivially, all arcwise connected metrizable spaces. There are also metrizable continua which are not arcwise connected and yet do have enough arcs (for instance that in Figure~\ref{fig:sine-curves} but also many indecomposable continua, including the solenoids and the Knaster type continua). Among the $\mathcal{AF}$-spaces one can find also spaces which do not have enough arcs (for instance the topologist's sine curve) and there are even disconnected $\mathcal{AF}$-spaces, including all zero-dimensional separable metrizable spaces, in particular, the Cantor space.

The main part of the paper asks what happens when we make new spaces out of old ones using standard topological constructions, such as sums (Section 3), products (Section 4), or continuous images and quotients (Section 5). We show by example that the $\mathcal{AF}$-property is not preserved under any of these operations, even in the category of compact metrizable spaces, see Theorems~\ref{thm:AFsum}, \ref{thm:AFproduct}, and \ref{thm:AFquotient}. However, we get positive results if we start with enough arcs. We show that an arbitrary topological sum of metrizable spaces with enough arcs is an $\mathcal{AF}$-space (although it is disconnected and so it does not have enough arcs), see Proposition~\ref{P:union}. Further, if a metrizable space $X$ is an (arbitrary) product of spaces with enough arcs, then it also has enough arcs and so is an $\mathcal{AF}$-space, see Proposition~\ref{prop:product} and Corollary~\ref{cor:product}. Finally, if the continuous image of a space with enough arcs is metrizable, then it is a space with enough arcs and so is an $\mathcal{AF}$-space, see Proposition~\ref{prop:quotient} and~Corollary~\ref{cor:quotient}. We remark that the subspace operation does not have a chance to behave well with respect to the $\mathcal{AF}$-property due to some standard embedding theorems (for instance every separable metrizable space of topological dimension $m$ can be embedded in $\mathbb R^{2m+1}$).

There are still some classes of continua not addressed by our methods. For example, we do not know if the pseudo-arc or pseudo-circle (or some other hereditarily indecomposable continua) are $\mathcal{AF}$-spaces; it would be interesting to give a complete topological characterization of $\alpha$-limit sets in these settings. We also did not systematically consider the  effect of passing to inverse limits. Another possible question is whether the hyperspace $2^X$ of an $\mathcal{AF}$-space $X$ is still an $\mathcal{AF}$-space, and in particular if all closed sets can be attained as $\alpha$-limit sets of ``hypermaps'' lifted from $X$.

\smallskip
 In Section~\ref{S:prelim} we present definitions, some observations and some examples. Then in Sections~\ref{S:sums}, ~\ref{S:prod} and ~\ref{S:quotients} we study whether the topological sum or product or quotient, respectively, of spaces with the $\mathcal {AF}$-property (or at least of spaces with enough arcs) has the $\mathcal {AF}$-property.

\section{Preliminaries and examples}\label{S:prelim}

\subsection{$\alpha$-limit sets and their basic properties}\label{SS:alpha}
If $X$ is a metrizable space and $f\colon X\to X$ is continuous, the $\alpha$-limit set of a point $x\in X$ is defined as
\begin{equation}\label{Eq:def alpha}
\alpha_f(x) := \bigcap_{n\geq 0} \overline{\bigcup_{k\geq n} f^{-k}(x)}. 
\end{equation}
Thus, $y\in \alpha_f(x)$ if and only if every neighbourhood of $y$ contains arbitrarily high preimages of $x$ (i.e., points $y_n$, $n=1,2, \dots$ with $f^{k_n}(y_n)=x$ for $k_n \nearrow +\infty$). Since we assumed metrizability (for simplicity, we make this assumption throughout the paper), each point $y\in X$ has a countable neighborhood basis and we have
\begin{equation}\label{Eq:alpha-1st count}
\alpha_f(x) = \{y\in X\colon (\exists y_n \to y) (\exists k_n\nearrow +\infty)(f^{k_n}(y_n) = x)\}.
\end{equation}
Let us summarize some basic facts on $\alpha$-limit sets.
\begin{enumerate}
	\item [(F1)] $\alpha_f(x)$ is closed.
	\item [(F2)] If $X$ is compact and $f$ is surjective, then $\alpha_f(x) \neq \emptyset$. 
	
	(Indeed, the intersection of a nested sequence of nonempty closed sets in a compact space is nonempty. Surjectivity without compactness is not sufficient -- to see this, consider $f(x)=x+1$ on $\mathbb Z$ or $\mathbb R$.)
	\item [(F3)] If $f$ is topologically exact, then $\alpha_f(x)=X$ for all $x\in X$.

    (The definition of topological exactness is that for each non-empty open set $U\subseteq X$ there is $n\geq0$ such that $f^n(U)=X$. Then~\eqref{Eq:def alpha} gives the result.)
	\item [(F4)] $f(\alpha_f(x)) \subseteq \alpha_f(x)$.
	
    (This inclusion may be strict even if we assume that $X$ is  compact or $f$ is surjective. Indeed, $X=\{0,1\}$ with the trivial metric is compact, the map $f$ such that $f(1)=f(0)=0$ is continuous, but not surjective, $\alpha_f(0) = \{0,1\}$ and $f(\alpha_f(0)) = \{0\}$. For the second counterexample with surjective $f$ and non-compact $X$, let $X$ be a subset of the plane defined as $X=\{\langle0,0\rangle\} \cup A \cup B \cup C$ where $A = \{\langle1/n, 0\rangle\colon n=1,2,\dots\} \cup \{\langle n,0\rangle\colon n=2,3,\dots\}$, $B=\{\langle1/n, 1/m\rangle\colon n=1,2, \dots \text{ and } m=n, n+1, \dots\}$ and $C=\{\langle n, 1 + 1/m\rangle\colon n=2,3,\dots \text{ and } m=1,2,\dots\}$. We put $f(\langle0,0\rangle)=\langle 0,0\rangle$ and define $f$ on $A$ by $f(\langle1/n, 0\rangle) = \langle1/(n+1),0\rangle$ for $n=1,2,\dots$ and $f(\langle n,0\rangle)=\langle n-1,0 \rangle$ for $n=2,3,\dots$. Further, we define $f$ on $B$ by $f(\langle1/n,1/n\rangle) =\langle 0,0 \rangle$, $n=1,2,\dots$ and $f(\langle1/n, 1/m\rangle) =\langle1/(n+1), 1/m \rangle$ for $n=1,2,\dots$ and $m=n+1, n+2, \dots$. Finally, to define $f$ on $C$, let $f(\langle2, 1+1/m\rangle) = \langle1,1/m \rangle$ for $m=1,2,\dots$ and $f(\langle n, 1+1/m\rangle) = \langle n-1,1+1/m\rangle$ for $n=3,4,\dots$ and $m=1,2,\dots$. Then $f\colon X\to X$ is continuous and surjective. One can see that $\langle 1,0\rangle \in \alpha_f(\langle 0,0\rangle)$ but the only preimage of $\langle 1,0\rangle$ is the point $\langle 2,0\rangle$ which does not belong to $\alpha_f(\langle 0,0\rangle)$. Therefore $f(\alpha_f( \langle 0,0\rangle)) \neq \alpha_f( \langle 0,0\rangle)$.) 
	
	\item [(F5)] If $X$ is compact and $f$ is surjective, then $f(\alpha_f(x)) = \alpha_f(x)$.
	
    (Indeed if $y\in\alpha_f(x)$ is the limit of points $y_n\to y$ with $f^{k_n}(y_n)=x$, $k_n \nearrow \infty$, then surjectivity guarantees preimages $z_n$ with $f(z_n)=y_n$ and compactness guarantees a convergent subsequence $z_{n_i} \to z \in X$. Clearly $f(z)=y$ and $z$ is in $\alpha_f(x)$).
	\item [(F6)] If $\alpha_f(x)$ has nonempty interior, then $x\in \alpha_f(x)$.
	
	(Let $c$ be an interior point of $\alpha_f(x)$. Then $\alpha_f(x)$, being a neighborhood of the point $c\in \alpha_f(x)$, contains arbitrarily high order preimages of $x$. However, by (F4), $\alpha_f(x)$ is invariant and so it contains also the point $x$.)
\end{enumerate}

\medskip

If $X$ is a metrizable space and $A \subseteq X$, we say that \emph{$A$ is an $\alpha$-limit set in $X$} provided that there is a continuous map $f\colon X\to X$ and a point $x\in X$ such that $A=\alpha_f(x)$. Denote by $\mathcal{A}(X)$ the family of $\alpha$-limit sets in $X$.
If $X$ is a singleton, then obviously $\mathcal A(X) = \{X\}$. 

The following is a trivial but useful observation for a dynamical system given by a space $X$ and a map $f\colon X\to X$. If, for a closed set $A\subseteq X$, we have $f(A)=\{a\} \subseteq A$ and $f(X\setminus A) \subseteq X\setminus A$, then $\alpha_f(a)=A$ (note that $a$ is a fixed point of $f$).

\begin{lemma}\label{L:trivial alpha-limit sets}
Let $X$ be nondegenerate. 
\begin{enumerate}
	\item $X\in \mathcal A(X)$ and $\emptyset \in \mathcal A(X)$.
	\item More generally, every clopen subset of $X$ belongs to $\mathcal A(X)$.	
	\item [(3)] Every singleton in $X$ belongs to $\mathcal A(X)$.
\end{enumerate} 	
\end{lemma}

\begin{proof}
(1) Choose two points $a\neq b$ in $X$. Let $f(x)=a$ for all $x\in X$. Then $\alpha_f(a) = X$ and 
	$\alpha_f(b) = \emptyset$.
	
(2) In view of (1), assume that $E$ is a nonempty, proper, clopen subset of $X$ and choose $e\in E$. Let $g(x)=e$ for all $x\in E$ and $g(x)=x$ for all $x\in X\setminus E$. Then $g$ is continuous and $\alpha_g(e)=E$.

(3) Fix $x_0\in X$ and let $h$ be the identity on $X$. Since every singleton in a metrizable space is a closed set,~\eqref{Eq:def alpha} gives that $\alpha_h(x_0)=\{x_0\}$.
\end{proof}

It is possible that only the sets listed in Lemma~\ref{L:trivial alpha-limit sets} belong to $\mathcal A(X)$.

\begin{example}[Rigid spaces]\label{Ex:rigid}
Recall that a nondegenerate space $X$ is \emph{rigid for continuous maps}, or just \emph{rigid}, if the only continuous maps $X\to X$ are the trivial ones, i.e., the identity and the constant maps. Clearly, a rigid space is connected. The first example of a rigid space was found by Cook~\cite{Cook} in the class of metrizable continua. (In fact Cook constructed in that paper what is now usually called a \emph{Cook continuum}. Roughly speaking, Cook continua are nondegenerate metrizable continua which are `everywhere' rigid. For a nontrivial application of Cook continua in dynamics see~\cite{SYZ}.) If $X$ is a rigid metrizable continuum, then obviously $\mathcal A(X)$ contains just the empty set, $X$, and all the singletons in $X$. Conversely, if $\mathcal A(X)$ contains just all these trivial sets, then $X$ is nondegenerate and connected; however, we do not know whether it has to be rigid.
\end{example}

\subsection{ $\mathcal A \mathcal F$-spaces and spaces with enough arcs}\label{SS:AF}
Denote by $\mathcal{F}(X)$ the family of closed sets in the metrizable space $X$. Recall that always $\mathcal A(X) \subseteq \mathcal F(X)$. 

\begin{definition}\label{D:AF}
In the case of equality $\mathcal{A}(X)=\mathcal{F}(X)$,  i.e. if the families of $\alpha$-limit sets and closed sets in $X$ coincide, we will say that the space $X$ has the \emph{$\mathcal A \mathcal F$-property} or that $X$ is an \emph{$\mathcal A \mathcal F$-space}. 
\end{definition}

As observed in~\cite{KMS},  all real compact intervals do have the $\mathcal A \mathcal F$-property. On the other hand, the singleton and, by Example~\ref{Ex:rigid}, the rigid metrizable continua do not have this property.
\begin{lemma}\label{lm:arc}
	Let $(X,d)$ be a metric space and $A\subseteq X$ a closed set such that there is an arc in $X$ whose one endpoint is in $A$ and the other one is in $X\setminus A$. Then $A\in \mathcal A(X)$.
\end{lemma}
\begin{proof}
Let $J=[a,b]$ be an arc with $a \in A$ and $b \in X\setminus A$ and $\gamma \colon [0,1]\to J$ be a homeomorphism with $\gamma(0)=a$ and $\gamma(1)=b$. We may assume that $J\cap A=\{a\}$, otherwise we  replace $J$ by the arc $[a', b]$ with $a'$ being the last point of $J$ lying in $A$, i.e.,  $a'=\gamma(s)$ where $s=\max\{t\in[0,1]: \gamma(t)\in A\}$. 

Define $f:X\to X$ by $f(x)=\gamma \left (\frac{d(x,A)}{d(x,A)+1}\right)$, $x\in X$. The map $f$ is continuous. If $x\in X\setminus A$ then $d(x,A)>0$, whence $\frac{d(x,A)}{d(x,A)+1}\in (0,1)$ and then $f(x)\in J\setminus\{a,b\}$. If $x\in A$ then $f(x)= \gamma(0) = a$.  Thus $f(A)=\{a\} \subseteq A$ and $f(X\setminus A)\subseteq X\setminus A$. It follows that $A=\alpha_f(a)$.
\end{proof}

The next definition is motivated by Lemma~\ref{lm:arc}. 
\begin{definition}\label{D:arcs}
Let $X$ be a nondegenerate metrizable space. If for every nonempty proper closed subset $A\subseteq X$ there is an arc whose one endpoint is in $A$ and the other one is in $X\setminus A$, then we will call $X$ a \emph{space with enough arcs}. 
\end{definition}
 Clearly, every nondegenerate arcwise connected metrizable space has enough arcs. The following proposition is a direct corollary of Lemmas~\ref{L:trivial alpha-limit sets}(1) and~\ref{lm:arc}.

\begin{proposition}\label{P:arcsAF}
	Every space with enough arcs has the $\mathcal A \mathcal F$-property.
\end{proposition}

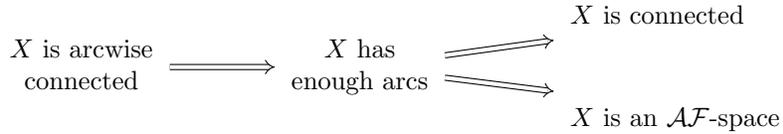
\begin{figure}[h!!]
\begin{tikzpicture}
\node[align=center,outer sep=3pt] (A) at (-4,0) {$X$ is arcwise\\ connected\vphantom{g}};
\node[align=center,outer sep=3pt] (B) at (-.3,0) {$X$ has\\ enough arcs};
\node[align=left, text width=3cm,outer sep=3pt] (C) at (4,0.7) {$X$ is connected};
\node[align=left, text width=3cm,outer sep=3pt] (D) at (4,-0.7) {$X$ is an $\mathcal{AF}$-space};
\draw[-implies,double equal sign distance] (B) -- (C.south west);
\draw[-implies,double equal sign distance] (B) -- (D.north west);
\draw[-implies,double equal sign distance] (A) -- (B);
\end{tikzpicture}
\caption{Implications in the class of nondegenerate metrizable spaces}\label{F:imps}
\end{figure}

Observe also that spaces with enough arcs are trivially connected. Then we arrive at the summary diagram shown in Figure~\ref{F:imps}.
We will see in the sequel that none of the reverse implications hold.

\begin{remark}\label{R:generalization}
In the present paper we restrict ourselves to metrizable spaces. The most important reason is that our main results, Theorems~\ref{thm:AFsum}, \ref{thm:AFproduct} and \ref{thm:AFquotient} are proved by finding examples and they all are metrizable. If in spite of this one wishes to consider the notions of an $\mathcal A \mathcal F$-space and of a space with enough arcs also for non-metrizable spaces, we wish to emphasize the following.  
		
\begin{itemize}
\item Definition~\eqref{Eq:def alpha} can be used in any topological space, but \eqref{Eq:alpha-1st count} works only in first countable spaces.
\item Facts (F1)-(F4) and (F6) hold in all topological spaces, (F5) in first countable sequentially compact spaces.
\item Lemma~\ref{L:trivial alpha-limit sets} works in any space.
\item Though the proof of Lemma~\ref{lm:arc} uses a metric, the only thing we really need is the existence of a continuous function $g: X\to [0,1]$ which vanishes precisely on~$A$; then put $f=\gamma \circ g$. If $X$ is normal, such a continuous function $g$ exists if and only if $A$ is a closed $G_{\delta}$ set in $X$, see~\cite[§14 VI, Theorems 1 and 2]{K1}, cf.~\cite[Exercise 4, p. 213]{Mu}. This shows that Lemma~\ref{lm:arc} can be extended to normal spaces in which all closed sets are $G_{\delta}$. Recall that such spaces are called \emph{perfectly normal}.
\item In Definitions~\ref{D:AF} and \ref{D:arcs} we need not assume metrizability. Then also Proposition~\ref{P:arcsAF} works without metrizability, provided it is assumed that the space is perfectly normal. The other two implications in Figure~\ref{F:imps} work for any space.
\end{itemize}
\end{remark}

\subsection{Examples of spaces with enough arcs (and so $\mathcal A \mathcal F$-spaces)}\label{SS:examples-enough}

 Now we give a criterion for a space to have enough arcs that does not require it to be arcwise connected.

\begin{proposition}\label{P:partition}
	Let $X$ be a nondegenerate metrizable space which can be written in the form $X=\bigcup_{\lambda \in \Lambda} X_{\lambda}$ with each $X_{\lambda}$ dense in $X$ and arcwise connected. Then $X$ has enough arcs.
\end{proposition}

\begin{proof}
 By Lemma~\ref{L:trivial alpha-limit sets}(1), $X\in \mathcal A(X)$ and $\emptyset \in \mathcal A(X)$. Let $A\subseteq X$ be a nonempty proper closed set. It intersects a set $X_{\lambda_0}$ for some $\lambda_0 \in \Lambda$ but it does not contain the whole $X_{\lambda_0}$  (otherwise we would have $A= \overline{A} \supseteq \overline {X_{\lambda_0}} = X$). Since $X_{\lambda_0}$ is arcwise connected, there is an arc in $X_{\lambda_0}$ joining a point of $A$ with a point not belonging to $A$. Thus $X$ has enough arcs. 
\end{proof}

\begin{corollary}\label{C:MinSolKnas}
The following metrizable spaces have enough arcs:
\begin{itemize}
	\item [(1)] Every nondegenerate metrizable space $X$ admitting a minimal continuous flow (= a continuous action of $\mathbb R$ on $X$).
	\item [(2)] Every nondegenerate metrizable continuum $X$ with all composants arcwise connected. In particular, this is the case for solenoids and Knaster continua.
\end{itemize}
\end{corollary}

\begin{proof}
(1) The orbits of the continuous flow form a partition of $X$. In general there are three types of orbits: points, circles and continuous injective images of $\mathbb R$.	Since $X$ is nondegenerate and the flow is minimal, no orbit is a point. Thus, again using minimality, either $X$ is a circle (and so $X$ has enough arcs) or $X$ decomposes into dense orbits which are continuous injective images of $\mathbb R$ (and so $X$ has enough arcs by Proposition~\ref{P:partition}).

(2) The composants of a nondegenerate metrizable continuum $X$ are dense (and connected) and their union is $X$ (see \cite[§48 VI]{K2}). By assumption, they are arcwise connected. So Proposition~\ref{P:partition} applies.

Solenoids/Knaster continua are inverse limits of circles/arcs with open bonding maps different from homeomorphisms (if all the bonding maps are the standard surjective tent maps with two linear pieces, we obtain the Knaster bucket handle continuum~\cite[Example 1, p. 204]{K2}). These nondegenerate metrizable continua are indecomposable (hence not arcwise connected) with composants coinciding with the arcwise components. (Each composant of a solenoid is an injective continuous image of the real line. The composants of a Knaster continuum are, with one or two exceptions, also of that kind. The exceptional composants are injective continuous images of a closed real half-line.)
\end{proof}

\begin{figure}[htb!!]
		\includegraphics[width=5cm]{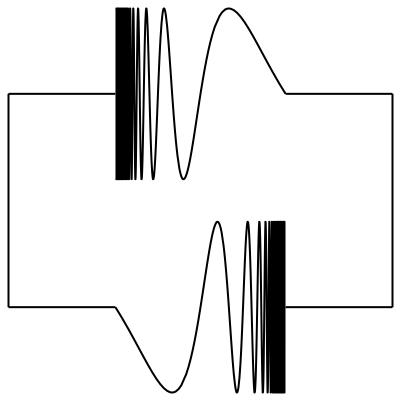}
		\caption{A space with enough arcs which is not arcwise connected}\label{fig:sine-curves}
	\end{figure}

We can now see that there are metrizable continua with enough arcs which are not arcwise connected. For instance \emph{indecomposable} continua (e.g. the solenoids and the Knaster continua) satisfying the assumptions of Corollary~\ref{C:MinSolKnas}(2) have these properties (recall that nondegenerate indecomposable continua are not arcwise connected~\cite[Exercise 11.54]{N}). There are also simpler examples, for instance we may glue two topologist's sine curves together as shown in Figure~\ref{fig:sine-curves} to get a \emph{decomposable} continuum with enough arcs which is not arcwise connected.

Besides nondegenerate arcwise connected continua, solenoids and Knaster continua, there are many other nondegenerate metrizable continua  which have all composants arcwise connected, and so are $\mathcal A \mathcal F$-spaces by Corollary~\ref{C:MinSolKnas}(2). For instance, this is true for inverse limits of certain unimodal bonding maps \cite{Ingram} or for inverse limit spaces of tent maps with finite critical orbit \cite{St}. By \cite[Proposition 4.9]{AABC} inverse limit spaces generated by long-branched tent maps have only arcs and points as proper subcontinua and thus they are $\mathcal A \mathcal F$-spaces. It would be interesting to \emph{characterize} those bonding maps $[0,1] \to [0,1]$ for which the corresponding inverse limit space has all composants arcwise connected (recall that the only arcwise connected inverse limit space (with onto bonding maps) of arcs is an arc~\cite{N0}). 

\subsection{Examples of $\mathcal A \mathcal F$-spaces lacking enough arcs}\label{SS:examples-AF}

What do we know about the class of spaces with the $\mathcal A \mathcal F$-property? First of all, the singleton and the rigid metrizable continua do not belong to this class. By Proposition~\ref{P:arcsAF} this class contains all spaces with enough arcs. We are going to show that this class is still larger and contains some well-known spaces lacking enough arcs. 

\begin{proposition}\label{P:0-dim}
Every nondegenerate zero-dimensional separable metric space $(X, d)$ (in particular, the Cantor space) has the $\mathcal A \mathcal F$-property.
\end{proposition}

\begin{proof}
Let $A$ be a nonempty, proper, closed subset of $X$. If $A$ is clopen then, by  Lemma~\ref{L:trivial alpha-limit sets}(2), it is an $\alpha$-limit set in $X$. 
Now assume that $A$ is not clopen. Then there is $a\in A$ such that 
\begin{equation}\label{Eq:0dist}
\dist(a,X\setminus A)=0.
\end{equation}
(Otherwise every point from $A$ has an open neighbourhood disjoint from $X\setminus A$; thus, this neighbourhood is a subset of $A$. It follows that the closed set $A$ is open, hence clopen, contradicting our assumption.) By~\eqref{Eq:0dist}, there is an injective sequence of points $b_i\in X\setminus A$, $i=1,2,\dots$, with $b_i\to a$. Choose pairwise disjoint clopen neighbourhoods $U_i$ of $b_i$, $i=1,2,\dots$, disjoint from $A$, with $\diam U_i \to 0$ (hence, the sequence $U_i$ converges to the point $a$). Then $A^* = A\sqcup \bigsqcup_{i=1}^\infty U_i$ is a closed set (here $\sqcup$ denotes disjoint union). By~\cite[p. 277]{K1}, every open set in $X$ is the union of a sequence of disjoint clopen sets (with arbitrarily small diameters). In particular, $X\setminus A^* = \bigsqcup_{i=1}^\infty V_i$ where the pairwise disjoint sets $V_i$ are clopen (possibly empty). Denote $B_i = U_i \sqcup V_i$, $i=1,2,\dots$. Then the sets $B_i$ are clopen, pairwise disjoint, disjoint from $A$ and so we have 
\[
X = A \sqcup \bigsqcup_{i=1}^\infty  B_i, \quad a\in A, \,\, b_i\in B_i, \, i=1,2,\dots, \,\, \text{ with } \,\, b_i\to a.
\]
We define a map $f\colon X\to X$ by $f(x)=a$ for $x\in A$, and $f(x)=b_i$ for $x\in B_i$. The map $f$ is continuous at all points from the clopen sets $B_i$. To prove the continuity at $x\in A$, fix a sequence $x_i \to x$. We want to prove that $f(x_i)\to f(x) =a$. Obviously, it suffices to prove this in two special cases, namely when all the points $x_i$ are in $A$ and when all the points $x_i$ are in $X\setminus A = \bigsqcup_{i=1}^\infty  B_i$.  If $\{x_i\}_{i=1}^{\infty}\subseteq A$ then $f(x_i)=a=f(x)$ and we are done. If $\{x_i\}_{i=1}^{\infty}\subseteq X\setminus A$ then $x_i\in B_{n(i)}$, $i=1,2,\dots$. Here $n(i) \to \infty$ because the point $x\in A$, which is the limit of $x_i$, does not belong to any of the closed sets $B_i$. It follows that $f(x_i) = b_{n(i)} \to a = f(x)$ and we are done again.
Thus the map $f$ is continuous. Clearly, $A=\alpha_f(a)$ because $f(A)=\{a\} \subseteq A$ and $f(X\setminus A)\subseteq X\setminus A$.
\end{proof}

\begin{figure}[htb!!]
\includegraphics[height=4.5cm]{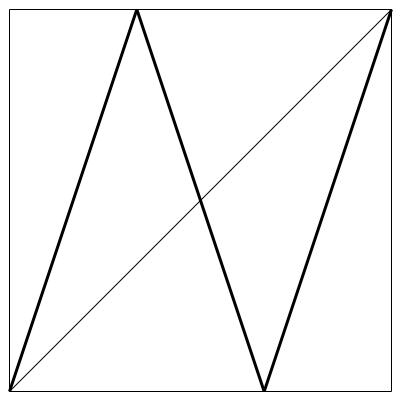}\hspace{1cm}
\includegraphics[height=4.5cm]{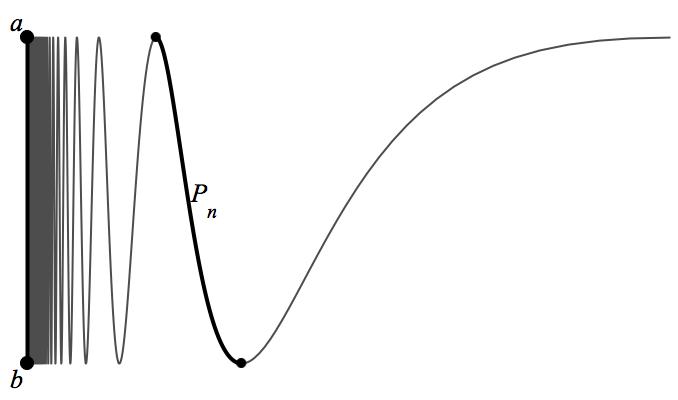}
\caption{The same map on each piece $P_n$ of the closed topologist's sine curve.}\label{fig:sine-curve}
\end{figure}

\begin{definition}\label{def:sine}[The sine curve map] Let $S=\{\langle x,\sin\frac{1}{x}\rangle: 0<x\leq\frac{2}{\pi}\}$ where $\langle a,b \rangle $ denotes the point in $\mathbb{R}^2$ with coordinates $a,b$. Then $\overline{S}\subseteq\mathbb{R}^2$ is the closed topologist's sine curve. We present a continuous map $g:\overline{S}\to\overline{S}$ which we will refer to henceforth as the sine curve map. Let $f:[-1,1]\to[-1,1]$ be the piecewise linear interval map with a graph determined by connecting the dots $\langle -1,-1\rangle ,\langle -\frac{1}{3}, 1\rangle , \langle \frac{1}{3},-1\rangle , \langle 1,1 \rangle$ and let $\pi_2:\bar{S}\to [-1,1]$ be the projection on the second coordinate, $\pi_2(\langle x,y\rangle)=y$, for every $\langle x,y\rangle\in \bar{S}$. We define $g$ on $\bar{S}$ piece-by-piece, see Figure~\ref{fig:sine-curve}. For every $n\in\mathbb{N}$, let $P_n$ be the unique arc in $S$ with endpoints $\langle \frac{2}{\pi(2n-1)},\sin (\frac{\pi(2n-1)}{2})\rangle$,  $\langle \frac{2}{\pi(2n+1)},\sin (\frac{\pi(2n+1)}{2})\rangle$ and define $g$ on $P_n$ by
\begin{equation}\label{eq:g}
 g|_{P_n}=(\pi_2|_{P_n})^{-1}\circ f \circ (\pi_2|_{P_n}).
\end{equation}
We use the same formula (\ref{eq:g}) for the arc $P$ in $\bar{S}$ with endpoints $a=\langle 0,1 \rangle$, $b=\langle 0,-1 \rangle $, which is the convergence continuum of $S$. The map $g$ is  continuous and well-defined since $\pi_2$ is a homeomorphism on every piece $P_n$ and the common endpoints of neighbouring pieces are fixed under $g$.
\end{definition}

\begin{lemma}
The closed topologist's sine curve has the $\mathcal A \mathcal F$-property but does not have enough arcs.
\end{lemma}

\begin{proof}
We keep the notation from Definition \ref{def:sine}. The arc $[a,b]$ forming the convergence continuum of $\bar{S}$ is not joined by any arc to its complement, so $\bar{S}$ is not a space with enough arcs. If $A$ is a non-empty proper closed subset of $\bar{S}$ other than the arc $[a,b]$ then there is an arc whose one endpoint is in $A$ and the other one is in $\bar{S}\setminus A$ and then we may apply Lemma \ref{lm:arc} to get that $A\in \mathcal A(\bar{S})$. To show that $[a,b]\in \mathcal A(\bar{S})$ and finish the proof we use the sine curve map $g$. Since the restriction $g|_{[a,b]}$ is topologically exact on the  arc $[a,b]$ we have $\alpha_g(b)\supseteq[a,b]$. Since $\bar{S}\setminus [a,b]$ is an open invariant set not containing $b$, it is disjoint from the $\alpha$-limit set of $b$. Therefore $\alpha_g(b)=[a,b]$.
\end{proof}

We give two more examples of spaces lacking enough arcs but with the $\mathcal A \mathcal F$-property. We will use them to show that $\mathcal A \mathcal F$-property is not preserved under topological sums (Theorem \ref{thm:AFsum}), products (Theorem \ref{thm:AFproduct}) or quotients (Theorem \ref{thm:AFquotient}).

\begin{figure}[htb!]
\includegraphics[height=4.5cm]{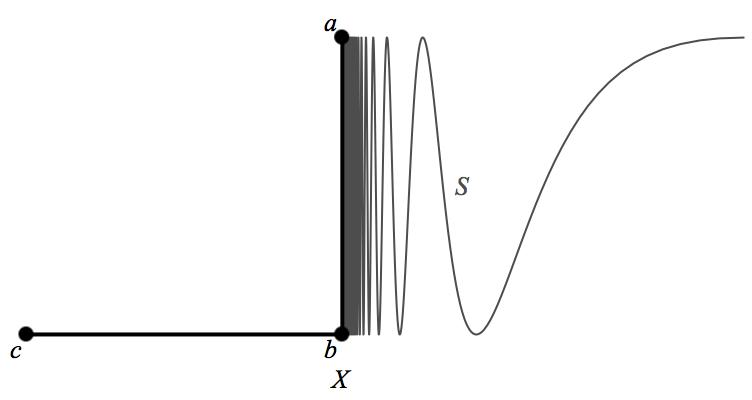}
\caption{The extended sine curve lacks enough arcs but has the $\mathcal A \mathcal F$-property.}\label{fig:X}
\end{figure}

\begin{definition}\label{def:X} [The extended sine curve]
 Let $a=\langle 0,1\rangle, b=\langle 0,-1\rangle, c=\langle -1,-1\rangle$ be points in $\mathbb{R}^2$. Let $X$ be formed by adjoining an arc $[b,c]$ to the convergence continuum $[a,b]$ of the topologist's sine curve $S=\{\langle x,\sin\frac{1}{x}\rangle: 0<x\leq\frac{2}{\pi}\}$, 
\[
X=\bar{S}\cup\{\langle x,-1 \rangle | -1\leq x \leq 0\}.
\]
We call $X$ the \emph{extended sine curve}, see Figure~\ref{fig:X}. Similarly as the topologist's sine curve, it clearly lacks enough arcs. Now extend the sine curve map $g:\bar{S}\to\bar{S}$ to a map $g:X\to X$ by putting $g([b,c])={b}$. We will refer to this map as the \emph{extended sine curve map}.  In Lemma~\ref{lem:XandY} we show that it is an $\mathcal A \mathcal F$-space.
\end{definition}

\begin{figure}[htb!]
\includegraphics[height=4.5cm]{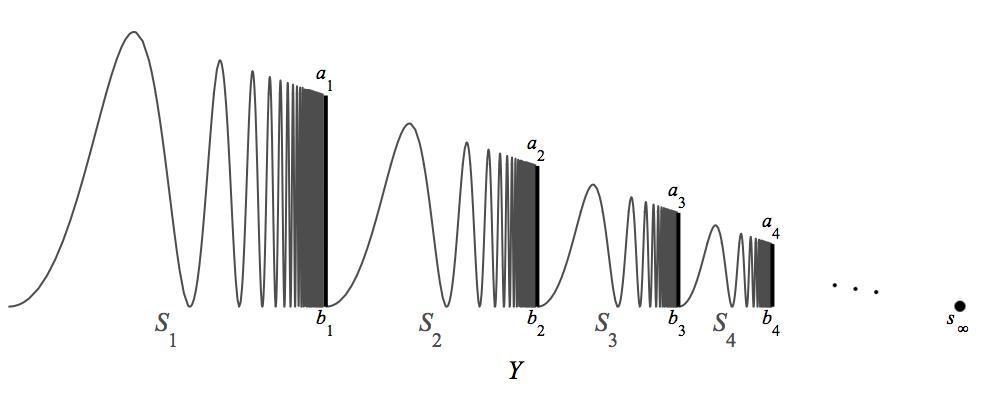}
\caption{A chain of sine curves lacks enough arcs but has the $\mathcal A \mathcal F$ property.}\label{fig:Y}
\end{figure}

\begin{definition}\label{def:Y} [The chain of sine curves]
Consider a topological space $Y$ which is the union of a sequence of closed topologist's sine curves $\bar{S_i}$ and a limit point $s_i$ not belonging to any of them and to which they converge. Suppose furthermore that the sets $\bar{S_i}$ are pairwise disjoint except for intersections $\bar{S}_i \cap \bar{S}_{i+1} = \{b_i\}$, where $\bar{S}_i = S_i \sqcup [a_i,b_i]$ and
\begin{enumerate}
\item $b_i$ is one of the endpoints of the arc $[a_i,b_i]$ forming the convergence continuum of $\bar{S_i}$.
\item $b_i$ is the unique endpoint of $S_{i+1}$, where we write $S_i = \bar{S_i} \setminus [a_i,b_i]$ for the open sine curves.
\end{enumerate}
The space $Y$ can be embedded in the plane $\mathbb{R}^2$, see Figure~\ref{fig:Y}. Note that, since  $s_{\infty}$ is a singleton, the space $Y$ is uniquely defined up to homeomorphism (in particular, though in Figure~\ref{fig:Y} the endpoints $b_i$ of $S_{i+1}$ are chosen to be the `bottom' endpoints of the convergence continua, for some and possibly for infinitely many $i$'s they could be the `top' endpoints); this fact will be used in the very end of the paper. We will refer to $Y$ as a \emph{chain of sine curves}. Clearly, $Y$ lacks enough arcs. By Lemma~\ref{lem:XandY}, it is an $\mathcal A \mathcal F$-space.
\end{definition}

\begin{lemma}\label{lem:XandY}
The extended sine curve and any chain of sine curves have the $\mathcal{AF}$-property.
\end{lemma}
\begin{proof}
We continue to use the notation from Definitions \ref{def:X} and \ref{def:Y}. So $X$ denotes the extended sine curve and $Y$ a chain of sine curves.

The only  nonempty proper closed subset of $X$ which is not joined in $X$ to its complement by an arc is the arc $[a,c]$, so in light of Lemma \ref{lm:arc} it suffices to show that $[a,c]\in \mathcal{A}(X)$. Let $g:X\to X$ be the extended sine curve map from Definition \ref{def:X}. Since the restriction of $g$ to the subspace $[a,b]$ is topologically exact we have $[a,b]\subseteq \alpha_{g}(b)$, and since the arc $[b,c]$ is collapsed to the fixed point $b$ we have $[b,c]\subseteq \alpha_{g}(b)$ as well. Moreover, $S$ is open, disjoint from $[a,c]$ and invariant for $g$, so $\alpha_{g}(b)\cap S=\emptyset$. Together this shows that $[a,c]= \alpha_{g}(b)$, completing the proof that $X$ has the $\mathcal{AF}$-property.

Any nonempty proper closed subset of $Y$ not joined by an arc to its complement is of the form $A_n:=[a_n,b_n]\cup \bigcup_{i>n}\bar{S_i}\cup\{s_{\infty}\}$, for some $n\in\mathbb{N}$, or is the singleton $\{s_{\infty}\}$. In light of Lemma \ref{lm:arc} it suffices to show that each of these sets belongs to $\mathcal{A}(Y)$. Singletons are always $\alpha$-limit sets by Lemma~\ref{L:trivial alpha-limit sets}. Now fix $n\in\mathbb{N}$. Let $f:Y\to Y$ be a map such that $f|_{\bar{S}_n}$ is conjugate to the sine curve map. Notice that $f(b_n)=b_n$. Put $f(\bigcup_{i>n}\bar{S_i}\cup \{s_{\infty}\})=b_n$. If $n\geq 2$, notice furthermore that $f(b_{n-1})=b_{n-1}$, so putting $f(\bigcup_{i<n}\bar{S_i})=b_{n-1}$ we get a continuous self-map on $Y$. We claim that $\alpha_f(b_n)=A_n$. Since $Y\setminus A_n$ is open, invariant and does not contain $b_n$, it must be disjoint from $\alpha_f(b_n)$. Since $f|_{\bar{S}_n}$ is the sine curve map, we get $[a_n,b_n]\subseteq \alpha_f(b_n)$. Finally since the rest of $A_n$ is collapsed to the fixed point $b_n$, we get  $\alpha_f(b_n)=A_n$. This completes the proof that $Y$ has the $\mathcal{AF}$-property.
\end{proof}

\section{Topological sums}\label{S:sums}
 If $\{X_i: i\in \Lambda\}$ is a family of any topological spaces then $\bigoplus_{i\in\Lambda} X_i$ denotes their topological sum. When working with it, we may without loss of generality assume that the underlying sets $X_i, i\in \Lambda$ are pairwise disjoint. We will use the symbol $\sqcup$ for disjoint union of sets.  Recall also that the topological sum of any family of metrizable spaces is metrizable.

If $X_1$ and $X_2$ are spaces with enough arcs, then their topological sum $X_1\oplus X_2$ does not have enough arcs (for $A=X_1$ there is no arc in $X$ joining a point from $A$ with a point from $X\setminus A$). However, we have the following observation. 

\begin{proposition}\label{P:union} 
Let $X_i$, $i\in \Lambda$, be spaces with enough arcs. Then their topological sum $X=\bigoplus_{i\in\Lambda} X_i$ has the $\mathcal A \mathcal F$-property. 
\end{proposition}

\begin{proof}
Let $A\subseteq X$ be a nonempty, proper, closed set. If there is $i\in\Lambda$ such that $\emptyset\neq A\cap X_i\subsetneq X_i$ then we apply Lemma~\ref{lm:arc} to get $A\in \mathcal A(X)$. Otherwise $A=\bigsqcup_{i\in\Lambda'} X_i$ where $\Lambda'$ is a nonempty, proper subset of $\Lambda$. Then $A$ is a clopen subset of $X$ and so $A\in \mathcal A(X)$ by Lemma~\ref{L:trivial alpha-limit sets}(2).
\end{proof}

Unfortunately, the $\mathcal{AF}$-property itself is not preserved under the topological sums, as shown by the next theorem.

\begin{theorem}\label{thm:AFsum}
There exists a pair of connected compact metrizable spaces each with the $\mathcal{AF}$-property such that their topological sum does not have the $\mathcal{AF}$-property.
\end{theorem}
\begin{proof}
Let $X$ be the extended sine curve and $Y$ be any chain of sine curves from Definitions~\ref{def:X} and~\ref{def:Y}. We continue to use the same notation as in those definitions. By Lemma \ref{lem:XandY}, $X$ and $Y$ have the $\mathcal{AF}$-property. We will show that $X\oplus Y$ does not have the $\mathcal{AF}$-property since $[a,c]\cup\{s_{\infty}\}\notin \mathcal{A}(X\oplus Y)$. 
\begin{claim}\label{claim:sum} $[a,c]\cup\{s_{\infty}\}\neq \alpha_f(x)$ for any pair $(f,x)$ where $x\in X\oplus Y$ and $f:X\oplus Y\to X\oplus Y$ is a continuous map.\end{claim}
Assume to the contrary that $[a,c]\cup\{s_{\infty}\}=\alpha_f(x)$, for some continuous map $f:X\oplus Y\to X\oplus Y$ and some point $x\in X\oplus Y$. Observe that $\alpha_f(x)$ is a neighborhood of $c$, so by fact (F6) from the beginning of Section~\ref{S:prelim} we get
\begin{equation}\label{E1}
x\in[a,c]\cup\{s_\infty\}.
\end{equation}
Since $\alpha_f(x)$ is invariant and connected components of an invariant set are always mapped into connected components, the arc $[a,c]$ is either invariant or is collapsed to the point $s_{\infty}$. We continue the proof in cases.\\
{\bf Case 1} $f([a,c])=s_{\infty}$\\
 The set $S$ from Definition~\ref{def:sine} is arcwise connected and with zero distance from the arc $[a,c]$, therefore $f(S)$ is an arcwise connected set with zero distance from the point $s_{\infty}$. The only such set is the singleton $\{s_{\infty}\}$. Therefore the whole set $X$ is collapsed to a single point in the first iterate. But $x$ has arbitrarily high order preimages in $X$ (because $X$ is a neighbourhood of a point from $\alpha_f(x)$), so the whole set $X$ must be mapped to $x$ infinitely often. Then $\alpha_f(x) \supset X$, a contradiction.\\
{\bf Case 2} $f([a,c])\subseteq [a,c]$ and $f(s_{\infty})=s_{\infty}$\\
Since $f(X),f(Y)$ are connected, we get $f(X)\subseteq X$ and $f(Y)\subseteq Y$. Thus $X,Y$ are invariant clopen sets and $x$ can have preimages in only one of them. Therefore $\alpha_f(x)$ is disjoint from one of the two sets $X,Y$. This is a contradiction.\\
{\bf Case 3} $f([a,c])\subseteq [a,c]$ and $f(s_{\infty})\neq s_{\infty}$\\
Since $\alpha_f(x)$ is invariant we have $f(s_{\infty})\in [a,c]$. Taking $W$ to be a small enough neighborhood of $f(s_\infty)$, the connected component of $W$ containing $f(s_\infty)$ is necessarily contained in $[a,c]$. Fix $N$ large enough that the set $V=\bigcup_{n\geq N}\bar{S}_n\cup\{{s_{\infty}}\}$ is mapped into $W$. Since $f(V)$ is connected and contains $f(s_{\infty})$, it follows that 
\begin{equation}\label{E2}
f(V)\subseteq [a,c].
\end{equation}
We will derive a contradiction by showing that $\alpha_f(x)$ contains other points besides $s_\infty$ in $V$. First observe that $x$ has arbitrarily high order preimages converging to $s_\infty$ in $V$, so if any iterated image of $V$ degenerates to a singleton, then $V$ maps to $x$ infinitely often giving $V\subset\alpha_f(x)$, a contradiction. Since $f([a,c])\subseteq [a,c]$, we conclude in particular that $f^2(V)$ is a non-degenerate connected subset of the arc $[a,c]$, i.e.\ a subarc. Then there must be $n\geq N$ such that $f^2(\bar{S}_n)$ is a non-degenerate subarc in $[a,c]$ as well, and in particular has non-empty interior relative to $[a,c]$. Choose $z$ from the relative interior and choose $w\in f(\bar{S}_n)$ so that $f(w)=z$. By~\eqref{E2}, $w \in f(\bar{S}_n) \subseteq f(V) \subseteq [a,c] \subseteq \alpha_f(x)$ and so there are points $w_i$ in $X$ and times $t_i\geq1$ such that
\begin{equation}\label{E3}
w_i \to w \text{ and } t_i \to \infty \text{ and } f^{t_i}(w_i) = x \text{ for all } i.
\end{equation}
We do not know if the points $w_i$ belong to $[a,c]$, some of them may belong to $S$. However, letting $z_i = f(w_i)$ we have
\begin{equation}\label{E4}
z_i \to z \text{ and } f^{t_i-1}(z_i) = x \text{ for all } i.
\end{equation}
Since $f(X)$ is connected and $[a,c]$ is invariant, we know that $X$ is invariant. Since $f(S)$ is arcwise connected it is either contained in $[a,c]$ or in $S$. In the first case the points $z_i$ clearly belong to $[a,c]$. In the second case $S$ is invariant, so by~\eqref{E1} and~\eqref{E3} it cannot contain any of the points $w_i$. Therefore the points $w_i$ and, since $f([a,c])\subseteq [a,c]$, also their images $z_i$ belong to $[a,c]$. Now since $z_i \to z$ in $[a,c]$ we get $z_i \in f^2(\bar{S}_n)$ for all sufficiently large $i$. By~\eqref{E4} it follows that $x \in f^{t_i+1}(\bar{S}_n)$ for all sufficiently large $i$, i.e.\ $x$ has arbitrarily high order preimages in $\bar{S}_n$. Then by compactness $\bar{S}_n$ contains a point from $\alpha_f(x)$, a contradiction.
\end{proof}

\section{Products}\label{S:prod}

Let $X=\prod_{\lambda\in\Lambda} X_{\lambda}$ be a product of topological spaces. Let $z=(z_{\lambda})_{\lambda\in \Lambda}\in X$ and $\lambda_0\in\Lambda$. The {\emph{line through $z$ along the coordinate $\lambda_0$} is the set $L(z,\lambda_0)=\{x\in X: x_{\lambda}=z_{\lambda}\text{ for all }\lambda\neq\lambda_0\}$. We will say that a point $a\in X$ is obtained from $b\in X$ by \emph{moving along a line} if there is a coordinate $\lambda_{0}\in \Lambda$ such that $a\in L(b,\lambda_0)$.
\begin{lemma}\label{lem:line}
Let $X=\prod_{\lambda\in\Lambda} X_{\lambda}$ be a product of topological spaces and $A\subset X$ be a nonempty, proper and closed set. Then there is a point $z\in A$ and a coordinate $\lambda_0\in\Lambda$ such that the line $L(z,\lambda_0)$ through $z$ along the coordinate $\lambda_0$ is not a subset of $A$.
\end{lemma}
\begin{proof}
Color the set $A$ white and the set $X\setminus A$ black. We will show that there is a point $z\in X$ and a coordinate $\lambda_0\in\Lambda$ such that the line $L(z,\lambda_0)$ is not monochromatic. Fix a white point $a=(a_{\lambda})_{\lambda\in \Lambda}\in A$ and a black point $b=(b_{\lambda})_{\lambda\in \Lambda}\in X\setminus A$. Since $X\setminus A$ is open, we can find a black neighborhood of $b$ of the form $B=\prod_{\lambda\in F}U_{\lambda}\times \prod_{\lambda\in\Lambda\setminus F}X_{\lambda}$ where $F\subset \Lambda$ is a finite set and $U_{\lambda}\subseteq X_{\lambda}$ is an open set for every $\lambda \in F$. Consider a point $c\in X$ with
$$
c_{\lambda}=
\begin{cases}
b_{\lambda}, \text{ for }\lambda \in F,\\
a_{\lambda},\text{ for }\lambda \in \Lambda\setminus F.
\end{cases}
$$
But $c$ can be obtained from the white point $a$ in $\# F$ steps by moving along a line (in each step we change one of the coordinates $a_{\lambda}$ to $b_{\lambda}$ for $\lambda\in F$). If all these lines are monochromatic then $c$ is a white point. But $c\in B$ by the definition of $B$, so $c$ is a black point. This shows that at least one of the lines is bichromatic.
\end{proof}
\begin{proposition}\label{prop:product} 
Let $X_{\lambda}$, $\lambda\in \Lambda$, be spaces with enough arcs. If their product $X=\prod_{\lambda\in\Lambda} X_{\lambda}$ is metrizable then it is a space with enough arcs.
\end{proposition}
\begin{proof}
Let $A$ be a nonempty proper closed subset of $X$. We need to find an arc joining a point from $A$ with a point from $X\setminus A$. By Lemma \ref{lem:line} there is a point $z\in A$ and a coordinate $\lambda_0$ such that the line $L(z,\lambda_0)\not\subset A$. Since $A$ is a closed set in $X$, $A\cap L(z,\lambda_0)$ is closed in $L(z,\lambda_0)$ equipped with the subspace topology. Let $\pi:X\to X_{\lambda_0}$ be the projection. Then the restriction $\pi|_{L(z,\lambda_0)}:L(z,\lambda_0) \to X_{\lambda_0}$ is a homeomorphism. Therefore $B=\pi(L(z,\lambda_0)\cap A)$ is a closed set in  $X_{\lambda_0}$. $B$ is nonempty because it contains $\pi(z)$ and is proper because $L(z,\lambda_0)\not\subset A$. Since $X_{\lambda_0}$ has enough arcs, there is an arc $\phi:[0,1]\to X_{\lambda_0}$ such that $\phi(0)\in B$ and $\phi(1)\not\in B$. Then $(\pi|_{L(z,\lambda_0)})^{-1}\circ \phi:[0,1]\to L(z,\lambda_0) \subseteq X$ is an arc in $X$ joining a point from $L(z,\lambda_0)\cap A$ to $L(z,\lambda_0)\setminus A$, that is from $A$ to $X\setminus A$.
\end{proof}
\begin{corollary}\label{cor:product}
If a product of spaces with enough arcs is metrizable, then it has the $\mathcal{AF}$-property. In particular, a countable product of metrizable spaces with enough arcs has the $\mathcal{AF}$-property.
\end{corollary}

\begin{theorem}\label{thm:AFproduct}
There exists a pair of compact metrizable spaces each with the $\mathcal{AF}$-property such that their product does not have the $\mathcal{AF}$-property.
\end{theorem}
\begin{proof}
Let $X$ be the extended sine curve from Definition \ref{def:X} and $Y$ be any chain of sine curves from Definition \ref{def:Y}. We continue to use the notation from those definitions. By Lemma \ref{lem:XandY}, $X$ and $Y$ have the $\mathcal{AF}$-property. First we will show that $X\oplus \{p\}$ and $Y\oplus \{p\}$ have also the $\mathcal{AF}$-property where $p$ is an arbitrary point disjoint from the spaces $X,Y$. 
\begin{claim} 
$X\oplus \{p\}$ has the $\mathcal{AF}$-property
\end{claim} 
The only proper closed subsets of $X\oplus \{p\}$ which are not joined in $X\oplus \{p\}$ to their complements by an arc are $[a,c]$, $\{p\}$ and the union $[a,c]\cup \{p\}$. In light of Lemma \ref{lm:arc} it suffices to show that these three sets belong to $\mathcal{A}(X\oplus \{p\})$. We extend the extended sine curve map $g$ from Definition~\ref{def:X} to maps $g_1,g_2:X\oplus\{p\}\to X\oplus\{p\}$ by setting $g_1(p)=p$ and $g_2(p)=b$. Obviously $\alpha_{g_1}(b)=[a,c]$ and $\alpha_{g_1}(p)=\{p\}$. Since $p$ is mapped by $g_2$ to the fixed point $b$ we have $\alpha_{g_2}(b)=[a,c]\cup \{p\}$.
\begin{claim} 
$Y\oplus \{p\}$ has the $\mathcal{AF}$-property
\end{claim}
Any proper closed subset of  $Y\oplus \{p\}$ not joined by an arc to its complements is one of the following sets: $\{s_{\infty}\}$, $\{p\}$,  $\{s_{\infty}\}\cup\{p\}$, $A_n$, $A_n\cup\{p\}$, for some $n\in\mathbb{N}$, where $A_n:=[a_n,b_n]\cup \bigcup_{i>n}\bar{S_i}\cup\{s_{\infty}\}$. Clearly $\{s_{\infty}\}=\alpha_{\id}(s_{\infty})$ and $\{p\}=\alpha_{\id}(p)$, where $\id$ denotes the identity. To obtain $\{s_{\infty}\}\cup\{p\}\in \mathcal{A}(Y\oplus\{p\})$ we define a map $h:Y\oplus\{p\}\to Y\oplus\{p\}$ such that $h|_{Y}=\id$ and $h(p)=s_{\infty}$. Then $\alpha_h(s_{\infty})=\{s_{\infty}\}\cup\{p\}$. Fix $n\in\mathbb{N}$. We extend the map $f$ acting on $Y$ from the proof of Lemma \ref{lem:XandY} to maps $f_1,f_2:Y\oplus\{p\}\to Y\oplus\{p\}$ by setting $f_1(p)=p$ and $f_2(p)=b_n$. It is easy to see that $\alpha_{f_1}(b_n)=A_n$ and  $\alpha_{f_2}(b_n)=A_n\cup\{p\}$ since $\alpha_f(b_n)=A_n$ and $b_n$ is a fixed point.

\begin{claim} 
	$W: = (X\oplus \{p\})\times (Y\oplus\{p\})$ does not have the $\mathcal{AF}$-property
\end{claim}

First observe that we may equivalently write 
\[
W = \underbrace{(X\times\{p\})\oplus (\{p\}\times Y)}_{W^{-}}\oplus \underbrace{(X\times Y)\oplus \{\langle p,p\rangle\}}_{W^{+}}.
\]
Thus $W$ is the topological sum of four pairwise disjoint continua, the last one being a singleton. We will show that the set 
\begin{equation}\label{Eq:V}
V:=([a,c] \times \{p\}) \cup (\langle p, s_{\infty}\rangle )  \, \subseteq \, W^{-} 
\end{equation}
does not belong to $\mathcal{A}(W)$. Assume to the contrary that 
\[
V=\alpha_F(w)
\] 
for some continuous map $F\colon W\to W$ and some point $w\in W$. The set $V$ lying in $W^{-}$ intersects both connected components $X\times\{p\}$ and $\{p\}\times Y$ of $W^{-}$. Moreover, being an $\alpha$-limit set, it is invariant. Therefore the $F$-images of those two components of $W^{-}$ are continua intersecting $W^{-}$. Since the clopen sets $W^{-}$ and $W^{+}$ forming $W$ are disjoint, no continuum in $W$ intersects both. Hence $F(W^{-}) \subseteq W^{-}$. Since points  from the clopen set $W^{+}$ have no preimages in the disjoint clopen set $W^{-}$ containing $\alpha_F(w)$, it follows that 
\[
w \in W^{-}.
\]
There is a natural homeomorphism $h\colon W^{-} \to  X\oplus Y$ defined by $h(\langle x, p\rangle )=x$ for every $x\in X$ and $h(\langle p, y\rangle )=y$ for every $y\in Y$. Since the continuous map $F|_{W^{-}}$ acts on $W^{-}$, we may define the corresponding conjugate map on $X\oplus Y$. So, we pass to the continuous map $f\colon X\oplus Y\to X\oplus Y$ defined by $f :=h\circ F|_{W^{-}} \circ h^{-1}$. We have $V=\alpha_F(w) = \alpha_{F|_{W^{-}}}(w) \subseteq W^{-}$ and so 
\[
h(V) = \alpha_f (h(w))
\]
is an $\alpha$-limit set in $X\oplus Y$. Using ~\eqref{Eq:V} and the definition of $h$ we get that $h(V) = [a,c]\cup\{s_{\infty}\}$. However, by Claim \ref{claim:sum} from the proof of Theorem \ref{thm:AFsum} this set is not an $\alpha$-limit set in $X\oplus Y$, a contradiction. 
\end{proof}

\section{Continuous images and quotient spaces}\label{S:quotients}

\begin{proposition}\label{prop:quotient} 
Let $X$ be a space with enough arcs and $f:X\to Y$ be a surjective continuous map. If the image space $Y$ is metrizable then  $Y$ has enough arcs.
\end{proposition}

\begin{proof}
Let $B\subset Y$ be a nonempty proper closed subset. Then the set $A:=f^{-1}(B)\subseteq X$ is a nonempty proper closed set and by assumption there is an arc $J_{A}$ joining a point $a\in A$ with a point $a'$ from the complement $X\setminus A$. The arc $J_A$ is given by a homeomorphism $h_A:[0,1]\to J_A$.  Thus $h_B=f\circ h_A$ is a continuous map from $[0,1]$ to $Y$ and $J_B=h_B([0,1])$ is a path in $Y$. Clearly $f(a)\in J_B\cap B$ and $f(a')\in J_B\cap (Y\setminus B)$. Since $J_B$ is a path-connected metric space there is an arc joining $f(a)\in B$ with $f(a')\in Y\setminus B$.
\end{proof}
\begin{corollary}\label{cor:quotient}
Let $X$ be a space with enough arcs. If $f$ is a continuous map such that $f(X)$ is metrizable, then $f(X)$ has the $\mathcal{AF}$-property. In particular, every metrizable quotient space $X/_{\sim}$ has  the $\mathcal{AF}$-property.
\end{corollary}

\begin{remark}\label{R:generalization2}
If one wishes to consider also non-metrizable spaces (cf. Remark~\ref{R:generalization}), here is a place to stop for a moment. In Proposition~\ref{prop:quotient} it is sufficient to assume that $Y$ is Hausdorff, since then $J_B$ is a Hausdorff path and it is well known that every Hausdorff path from $x$ to $y$ contains an arc from $x$ to~$y$. Then in Corollary~\ref{cor:quotient} it is sufficient to assume that $f(X)$ is Hausdorff and perfectly normal. In particular, every Hausdorff perfectly normal quotient space of a space with enough arcs has  the $\mathcal{AF}$-property.
\end{remark}

\begin{theorem}\label{thm:AFquotient}
There is a compact metrizable space $Z$ with the $\mathcal{AF}$-property and an equivalence relation $\sim$ on $Z$ such that the quotient space $Z/_{\sim}$ is a compact metrizable space which does not have the $\mathcal{AF}$-property.
\end{theorem}
\begin{proof}
Let $X\subseteq \mathbb{R}^2$ be the extended sine curve. Recall that $a=\langle 0,1\rangle, b=\langle 0,-1\rangle, c=\langle -1,-1\rangle$ and $X$ is formed by adjoining an arc $[b,c]$ to the convergence continuum $[a,b]$ of the topologist's sine curve $S=\{\langle x,\sin\frac{1}{x}\rangle: 0<x\leq\frac{2}{\pi}\}$,
$$X=\bar{S}\cup\{\langle x,-1 \rangle | -1\leq x \leq 0\}.$$
For $n\in\mathbb{N}$, $P_n$ is the unique arc in $S$ with endpoints $\langle \frac{2}{\pi(2n-1)},\sin (\frac{\pi(2n-1)}{2})\rangle$,  $\langle \frac{2}{\pi(2n+1)},\sin (\frac{\pi(2n+1)}{2})\rangle$. We call these sets the ``pieces'' of the sine curve. 

Let $W\subseteq \mathbb{R}^2$ be defined as follows. For every $i\in\mathbb{N}$, let $S_i=M_i(S)$, $a_i=M_i(a)$, $b_i=M_i(b)$, where $M_i$ is the affine mapping $M_i(\langle x,y\rangle)=\langle -2-\frac{1}{2^i}-\frac{(\pi/2) x}{2^i},(-1)^{i+1} y\rangle$ and let $S_{\infty}=\{\langle -2, y\rangle : -1\leq y\leq 1\}$. Define  
\[
W=S_{\infty}\cup \bigcup_{i\geq 1}\bar{S}_i = S_\infty \cup \bigcup_{i\geq 1}M_i(\bar{S}).
\]

Finally, let $Z=W\cup X$ be the union in $\mathbb{R}^2$ of these two compact metrizable spaces, see Figure~\ref{fig:Z}. It wil be convenient to use, for $n\in \mathbb N$, the notation 
\[
A_n:=[a_n,b_n]\cup  \bigcup_{i>n}\bar{S_i}\cup S_{\infty}.
\] 
As usual, by $[a,c]$ we denote the unique arc in $Z$ with endpoints $a,c$. 

\begin{figure}
\includegraphics[width=.9\textwidth]{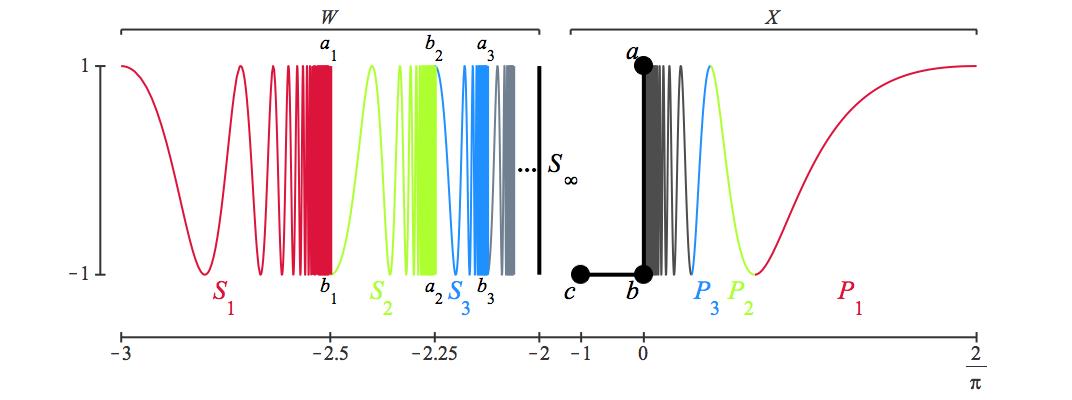}
\caption{Passing to quotients does not preserve the $\mathcal{AF}$-property (horizontal scaling is distorted to show more detail).}\label{fig:Z}
\end{figure}

We will show that $Z$ has the $\mathcal{AF}$-property. In light of Lemma \ref{lm:arc}, it is sufficient to show that any nonempty proper closed subset of $Z$ not joined by an arc to its complement belongs to $\mathcal{A}(Z)$. Moreover, by Lemma~\ref{L:trivial alpha-limit sets}(2), among such sets we need not consider the clopen sets $W$ and $X$. Thus we come to the following list of sets (with $n\in \mathbb N$) which are to be considered.  
	
\begin{itemize}
	\item (Connected subsets of $W$) \,\,  $A_n$, $S_{\infty}$,
	\item (Connected subsets of $X$) \,\,  $[a,c]$,	
	\item (Disconnected sets) \,\, $W \cup [a,c]$, $A_n \cup X$, $A_n \cup [a,c]$, $S_{\infty}\cup X$, $S_{\infty}\cup [a,c]$.
\end{itemize}
We are going to prove that each of them equals $\alpha_f(z)$ for some continuous selfmap $f$ of $Z$ and some point $z\in Z$.

$\bullet$ $A_n \in \mathcal{A}(Z)$: The map $f$ is the identity on $X$ and on $W$ it is similar to the map from the proof of Lemma \ref{lem:XandY}. Thus, on $\bar{S}_n$ we put a copy of the sine curve map, we collapse $\bigcup_{i<n} \bar{S}_i$ to the fixed point $b_{n-1}$ (if $n\geq2$) and collapse $\bigcup_{i>n} \bar{S}_i \cup S_\infty$ to the fixed point $b_n$. Then $\alpha_f(b_n) = A_n$.

$\bullet$ $S_{\infty} \in \mathcal{A}(Z)$: Let $f$ be the identity on $X$ and let $g$ be the sine curve map on $\bar{S}$. On $\bar{S_i}$ let $f$ be the composition $M_i \circ g \circ M_i^{-1}$ (a conjugate copy of the sine curve map). Note that all the points $b_i$ are fixed for a map defined in this way, so this definition is consistent. Continuity dictates that the map $f$ on $S_\infty$ is the full 3-horseshoe (note that the 3-horseshoe is conjugate to itself when we reverse the orientation of $[-1,1]$ via the ``flip'' $y\mapsto -y$). Then for any point $x\in S_\infty$ the alpha-limit set of $x$ contains $S_\infty$ by topological exactness and is disjoint from the rest of $Z$ since $Z\setminus S_\infty$ is an open invariant set.

$\bullet$ $[a,c] \in \mathcal{A}(Z)$: Let $f:Z \to Z$ be the extended sine curve map on $X$ and the identity on $W$. Clearly $\alpha_f(b)=[a,c]$.

$\bullet$ $W \cup [a,c] \in \mathcal{A}(Z)$: Let $f:Z \to Z$ be the extended sine curve map on $X$ and let $f(W) = \{b\}$. Then $\alpha_f(b)=W \cup [a,c]$.

$\bullet$ $A_n \cup X \in \mathcal{A}(Z)$: On $W$ let $f$ be the same as in the case $A_n\in \mathcal A(Z)$ and let $f(X) = \{b_n\}$. Then $\alpha_f(b_n) = A_n \cup X$.

$\bullet$ $A_n \cup [a,c] \in \mathcal{A}(Z)$: Define $f$ to be the extended sine curve map on $X$. On $W$ we define $f$ in pieces. On $\bar{S_n}$ we take $f$ to be the homeomorphism $M_n^{-1}:\bar{S}_n\to \bar{S}$, which is a homeomorphism sending the convergence continuum $[a_n,b_n]$ of $S_n$ onto $[a,b]$ by a horizontal translation, followed by the flip in the second coordinate $y\mapsto -y$ if $n$ is even, and sending the endpoint $b_{n-1}$ of $S_n$ to the endpoint of $S$. Further, if $n$ is odd let $f$ send $A_n$ onto $[a,b]$ by horizontal projection, preserving the second coordinate, and if $n$ is even let $f$ send $A_n$ onto $[a,b]$ by horizontal projection  followed by the flip in the second coordinate. If $n\geq2$ then let $f$ collapse the set $\bigcup_{i<n}\bar{S}_i$ to the point $f(b_{n-1})$, which was already defined above as the endpoint of $S$. We claim that $\alpha_f(b)=A_n \cup [a,c]$. Clearly the $\alpha$-limit set contains $[a,c]$ since the map $f$ is the extended sine curve map on $X$. By topological exactness there are preimages of $b$ of arbitrarily high order in every subarc of the arc $[a,b]$, and every such preimage has preimages at all points of $A_n$ with the same second coordinate (when $n$ is odd) or with flipped second coordinate (when $n$ is even). This shows that the $\alpha$-limit set contains $A_n\cup [a,c]$. The rest of the space $Z$ is the union $S\cup S_n\cup \bigcup_{i<n}\bar{S}_i$ and it is again an open invariant set not containing $b$, and therefore disjoint from $\alpha_f(b)$.

$\bullet$ $S_{\infty}\cup X \in \mathcal{A}(Z)$:  On $W$, let $f$ be defined in the same way as in the case $S_{\infty} \in \mathcal{A}(Z)$. Further, let $f$ send $X$ to a fixed point $w$ of $f$ in $S_{\infty}$. Then $\alpha_f(w)= S_{\infty}\cup X$.
 
$\bullet$ $S_{\infty}\cup[a,c] \in \mathcal{A}(Z)$: We define a map $f: Z\to Z$ such that $S_{\infty}\cup[a,c]=\alpha_f(b)$. For $\langle x,y \rangle \in \bar{S}_i$ let $f(\langle x,y \rangle)$ be the unique point in $P_i$ with the second coordinate $y$. Note that this is well-defined since the point $b_i$ at the intersection of $\bar{S}_i$ with $\bar{S}_{i+1}$ is mapped to the point of intersection of $P_i$ with $P_{i+1}$ (this is the reason why $M_i$ has the coefficient $(-1)^{i+1}$). Then continuity dictates that $f$ carries $S_{\infty}$ by translation onto the arc $[a,b]$ in $X$. For defining $f$ on $X$ we simply use the extended sine curve map. Then $\alpha_f(b)$ contains $[a,b]$ since $f|_{[a,b]}$ is a topologically exact interval map, it contains $[b,c]$ since this whole arc is collapsed to the fixed point $b$ and it contains $S_{\infty}$ since each preimage of $b$ in $[a,b]$ translates to a preimage in $S_{\infty}$ with the same second coordinate and with the order of the preimage increased by 1. The complement of $S_{\infty}\cup[a,c]$ in $Z$ is an open invariant set not containing $b$, hence it is disjoint from $\alpha_f(b)$ and $S_{\infty}\cup[a,c]=\alpha_f(b)$.

\smallskip

We have now checked that all closed subsets of $Z$ are alpha-limit sets, so $Z$ has the $\mathcal{AF}$-property. However, by the proof of Theorem \ref{thm:AFsum}, the topological sum $X\oplus Y$ does not have the $\mathcal{AF}$-property, where $X$ is the extended sine curve and $Y$ is any chain of sine curves. Let $\sim$ be the equivalence relation on $W$ which identifies all points of $S_\infty$ to a single point $s_\infty$ and does not make any other identifications. Then $W/_{\sim}$ satisfies Definition \ref{def:Y} and therefore it is a chain of sine curves. Simultaneously $X\oplus W/_{\sim}$ is a factor of $Z$ which completes the proof. 
\end{proof}

\addtocontents{toc}{\protect\setcounter{tocdepth}{0}}
\subsection*{Acknowledgements}
The authors thank an anonymous referee for suggestions how to improve the presentation of the paper.
The second author was supported by RVO, Czech Republic funding for I\v{C}47813059. The third author was supported by the Slovak Research and Development Agency under contract No.~APVV-15-0439 and by VEGA grant 1/0158/20.

\begin{table}[h]

\begin{tabular}[t]{b{1.5cm} m{13.5cm}}

\includegraphics [width=.09\textwidth]{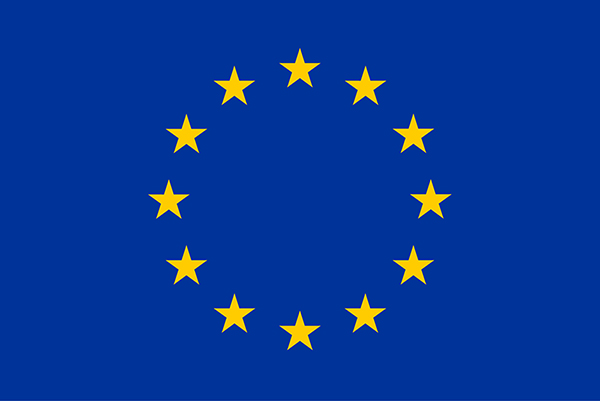} & 
This research is part of a project that has received funding from the European Union's Horizon 2020 research and innovation programme under the Marie Sk\l odowska-Curie grant agreement No 883748.

\end{tabular}
\end{table}

\end{document}